\documentclass[12pt,a4paper,dvips]{article}
\usepackage{inputenc, amssymb, amsmath, amsthm, srcltx}
\usepackage[english]{babel}
\usepackage{color}
\usepackage{graphicx}
\makeatletter
\thm@style{plain}

\def\th@plain{%
  \thm@headfont{\bfseries}%
  \itshape % body font
  \thm@notefont{\rm}%
}
\def\thm@indent{\hspace*{\parindent}}
\makeatother

\def\({\left(}
\def\){\right)}

\newcommand{\be}{\begin{equation}}
\newcommand{\ee}{\end{equation}}

\renewcommand{\leq}{\leqslant}

\let\epsilon\varepsilon
\let\phi\varphi
\let\le\leqslant
\let\ge\geqslant

\def\bydef{\ensuremath{\mathrel{{:}\,{=}}}}

\newtheorem{theorem}{Theorem}
\newtheorem{lemma}{Lemma}
\newcommand{\dokvo}{{\it Proof.} }
\newcommand\inte{\int\limits}

\newtheorem{conj}{Conjecture}

\newcommand{\beq}{\begin{equation}}
\newcommand{\eeq}{\end{equation}}

{Definition}
{Algorithm}

\begin{document}

\centerline{\bf\uppercase{Correlations of the Moebius and Liouville functions}}
\centerline{\bf\uppercase{and the twin prime conjecture}\footnote[1]{%
}}

\bigskip

\medskip

\centerline{\sc Sergei~Preobrazhenski\u i and Tatyana~Preobrazhenskaya}

\bigskip

\bigskip

\hbox to \textwidth{\hfil\parbox{0.9\textwidth}{%
\small {\sc Abstract.} In this note we describe
weight functions that exhibit a transitional behavior between weak and strong correlation
with the Liouville function.
We also describe a binary problem which may be considered as an interpolation between
Chowla's conjecture for two-point correlations of the M\"obius function
and the twin prime conjecture,
in view of recent parity breaking results of K.~Matom\"aki, M.~Radziwi{\l}{\l} and T. Tao.}\hfil}

\bigskip

\bigskip

%------------------------------------------------------------------

{\bf 1. Introduction.} The Chowla conjecture for the Liouville function $\lambda(n)$ asserts that
\[
\sum_{n\le N}\lambda^{e_1}(n+h_1)\dots\lambda^{e_k}(n+h_k)=o(N)
\]
as $N\to\infty$ for any fixed integers $h_1$, $\ldots$, $h_k$ with $h_i\ne h_j$ for $i\ne j$
and with at least one odd exponent $e_i$.
In~\cite{PreobTao15} Tao proved the following logarithmically averaged version of Chowla's conjecture for two-point correlations:
\begin{theorem}[T. Tao]
Let $a_1$, $a_2$ be natural numbers, and let $b_1$, $b_2$ be integers such that $a_1b_2-a_2b_1\neq0$.
Let $1\leq h(x)\leq x$ be a quantity depending on $x$ that goes to infinity as $x\rightarrow\infty$.
Then one has
\[
\sum_{x/h(x)<n\leq x}\frac{\lambda(a_1n+b_1)\lambda(a_2n+b_2)}{n}=o(\log h(x))
\]
as $x\rightarrow\infty$.
\end{theorem}

The twin prime conjecture is the assertion that
\[
H_1\bydef\liminf\limits_{n\to\infty}(p_{n+1}-p_n)=2.
\]
In 2004 Goldston, Pintz and Y{\i}ld{\i}r{\i}m~\cite{PreobGPY09} established $H_1\le16$
on the Elliott--Halberstam conjecture~\cite{PreobEH70}.
The breakthrough paper of Zhang~\cite{PreobZh13} shows that $H_1\le70\,000\,000$ unconditionally.
Subsequent improvements~\cite{PreobPm14} lowered this bound to $H_1\le246$.
From the generalized Elliott--Halberstam hypothesis (see~Conjecture~\ref{PreobGEH}) Polymath8b project
has managed to reach $H_1\le6$. The paper~\cite[Section 8]{PreobPm14} also provides a heuristic argument that
the parity problem prohibits a sieve-theoretic proof of the twin prime conjecture.

Now we state a version of the generalized Elliott--Halberstam conjecture from~\cite{PreobPm14}.
The \emph{Dirichlet convolution} $\alpha\star\beta\colon\mathbb{N}\to\mathbb{C}$
of two arithmetic functions $\alpha,\beta\colon\mathbb{N}\to\mathbb{C}$ is defined as
\[
\alpha\star\beta(n)\bydef\sum_{d|n}\alpha(d)
\beta\left(\frac{n}{d}\right)
=\sum_{ab=n}{\alpha(a)\beta(b)}.
\]
For any function $\alpha\colon\mathbb{N}\to\mathbb{C}$ with finite support
(that is, $\alpha$ is non-zero only on a finite set) and any primitive congruence class $a\ (q)$,
we define the (signed) \emph{discrepancy} $\Delta(\alpha;a\ (q))$ to be the quantity
\[
\Delta(\alpha;a\ (q))\bydef\sum_{n=a\ (q)}\alpha(n)
-\frac{1}{\phi(q)}\sum_{(n,q)=1}\alpha(n).
\]
\begin{conj}[GEH $(D)$]
\label{PreobGEH}
Let $D=x^{1-\varepsilon(x)}$ for some function $\varepsilon(x)>0$ and $x/D\leq N,M\leq D$, be such that $NM \sim x$,
and let $\alpha, \beta$ be coefficient sequences at scale $N,M$. Then
\[
\sum_{q \lessapprox D} \sup_{a \in ({\mathbb Z}/q{\mathbb Z})^{\times}} |\Delta(\alpha\star\beta; a\ (q))| \ll x \log^{-A} x
\]
for any fixed $A>0$.
\end{conj}
The generalized Elliott--Halberstam conjecture (GEH) will refer to the assertion that GEH $(D)$ holds for all
$D=x^{1-\varepsilon(x)}$ with $\frac1{(\log x)^{1-\delta}}\le\varepsilon(x)<1$ and some fixed $0<\delta<1$.
We remark that from the results of~\cite{PreobFGHM91} it follows that Conjecture~\ref{PreobGEH}
fails with $\varepsilon(x)\le\frac{(A-\delta)(\log\log x)^2}{(\log x)(\log\log\log x)}$.

Let $\lambda(n)$ be the Liouville function, $\mu(n)$ be the M\"obius function,
$\tau_r(n)$ be the $r$-th divisor function for real $r$, that is the multiplicative function given
for all primes $p$ and all $m\ge1$ by
\[
\tau_r(p^m)=\frac{\Gamma(r+m)}{\Gamma(r)m!}.
\]

For a positive integer $n$ and $y\ge2$
let $\omega{{+}}(n,y)$ be the number of prime divisors of $n$ which are
$\ge y$ and $\omega{{-}}(n,y)$ be the number of prime divisors of $n$ which are
$<y$. On squarefree integers $n$ define the multiplicative function $\tau^{{\pm}}_{\kappa_1,\kappa_2}(n,y)$ by
\[
\tau^{{\pm}}_{\kappa_1,\kappa_2}(n,y)\bydef\kappa_1^{\omega{{-}}(n,y)}\kappa_2^{\omega{{+}}(n,y)}.
\]

Let $1_{S}(n)$ be the indicator function of $S$.
Suppose that $(h_1,h_2)$ is a fixed admissible $2$-tuple.

Another assumption of this note is based on the relationship between
the function $\tau^{{\pm}}_{\kappa_1,\kappa_2}(n,y)$ with some fixed $0\le\kappa_1<\kappa_2$ and
the integers left unsieved when sieving with the primes less than
$y=\exp\left((\log x)^{\delta}\right)$ for some fixed $0<\delta<1$.
More precisely, we have the following relevant results.
\begin{theorem}[K. Alladi~\cite{PreobAll82}]
Define
\[
S(x,y)\bydef\{n\le x\ {:}\ \text{\textup{least prime divisor of $n$ is $\ge y$}}\},
\]
\[
\phi_t(x,y)\bydef\#\{n\in S(x,y)\ {:}\ \omega(n)-\log u<t\sqrt{\log u}\},
\]
where $u=\log x/\log y$. Then for $2\le y\le x$
\[
\left|\phi_t(x,y)/|S(x,y)|-\frac1{\sqrt{2\pi}}\inte_{-\infty}^te^{-v^2/2}\,dv\right|\ll\frac1{\sqrt{\log u}}.
\]
\end{theorem}
\begin{theorem}[A. Selberg] %~\cite{PreobSel54}
Let $\pi_k(x)$ be the number of positive integers less than $x$ which have exactly $k$
distinct prime divisors (not counting multiplicity). Then for $k=O(\log\log x)$
\[
\pi_k(x)\sim F\left(\frac{k-1}{\log\log x}\right)\frac x{\log x}\frac{(\log\log x)^{k-1}}{(k-1)!},
\]
where
\[
F(z)=\frac1{\Gamma(z+1)}\prod_p\left(1+\frac z{p-1}\right)\left(1-\frac1p\right)^z
\]
and for $0<\kappa\le1$
\[
\sum_{n\le x}\tau_{\kappa}(n)\sim c(\kappa)x(\log x)^{\kappa-1}.
\]
\end{theorem}
In view of this results, one may consider the weight $\tau_{\kappa}(n)$ for $0<\kappa<1$ as a kind
of preliminary sieving, since this weight is mostly concentrated on numbers
which have $\kappa\log\log x$ prime divisors.
Next we can make actual sieving in the weighted sum
\[
\sum_{N/2<n\leq N}\lambda(n+h_1)\lambda(n+h_2)\tau_{\kappa}(n+h_1)\tau_{\kappa}(n+h_2)
\]
with the primes less than $y=\exp\left((\log x)^{\delta}\right)$.
This construction seems to be more viable when $\kappa$ is close to $1$ and $\delta$ is close to $0$.

More generally, we can make the following
\begin{conj}[Generalized Chowla's conjecture]
\label{PreobGCC}
For any fixed $0\le\kappa_1<\kappa_2$ and
$y=\exp\left((\log x)^{\delta}\right)$ with some fixed $0<\delta<1$ we have
\[
\sum_{\frac N{2}\le n\le N}\mu(n+h_1)\mu(n+h_2)
\tau^{{\pm}}_{\kappa_1,\kappa_2}(n+h_1,y)\tau^{{\pm}}_{\kappa_1,\kappa_2}(n+h_2,y)
=o\left(
\sum_{\frac N{2}\le n\le N}
\tau^{{\pm}}_{\kappa_1,\kappa_2}(n+h_1,y)\tau^{{\pm}}_{\kappa_1,\kappa_2}(n+h_2,y)
\right).
\]
\end{conj}
In the case $\kappa_1=\kappa_2=1$ this is the original Chowla's conjecture
for two-point correlations of the M\"obius function.
And in the case of any $0\le\kappa_1<\kappa_2$ this conjecture essentially implies
the twin prime conjecture in view of the results of Bombieri (see~\cite[Chapter 16]{PreobFI10}).
Obviously, if we had the case $\kappa_1=\kappa_2=1$ then by continuity we could perturb
$\kappa_1$, $\kappa_2$, $y$ to get into the region $\kappa_1<\kappa_2$.

Next, define the weight $\alpha_{\kappa}(n)$ for $0\le\kappa\le1$ by
\[
\alpha_{\kappa}(n)
\bydef\sum_{d\mid n}\mu(d)\tau_{1+\kappa}\left(\frac nd\right)\left(\frac{\log n}{\log N}-\frac{\log d}{\log N}\right).
\]
The result of this note is
\begin{theorem}\label{Preoblambdacorrel}
We have
\[
\sum_{\frac N{2}\le n\le N}\lambda(n)\alpha_{1}(n)
=o\left(\sum_{\frac N{2}\le n\le N}\alpha_{1}(n)\right)
\]
but
\[
\sum_{\frac N{2}\le n\le N}\lambda(n)\alpha_{0}(n)
=-\left(\sum_{\frac N{2}\le n\le N}\alpha_{0}(n)\right)
\]
with the sums on the right-hand sides being nonzero for all sufficiently large values of $N$.
For two-point correlations, we have unconditionally
\[
\sum_{\frac N{h(N)}\le n\le N}\frac{\lambda(n+h_1)\lambda(n+h_2)\alpha_{1}(n+h_1)\alpha_{1}(n+h_2)}{n}
=o\left(\sum_{\frac N{h(N)}\le n\le N}\frac{\alpha_{1}(n+h_1)\alpha_{1}(n+h_2)}{n}\right)
\]
and assuming Conjectures~\textup{\ref{PreobGEH}--\ref{PreobGCC}}
\[
\sum_{\frac N{2}\le n\le N}\lambda(n+h_1)\lambda(n+h_2)\alpha_{0}(n+h_1)\alpha_{0}(n+h_2)
=\sum_{\frac N{2}\le n\le N}\alpha_{0}(n+h_1)\alpha_{0}(n+h_2)
\]
with the sum on the right-hand side being nonzero.
The latter claim implies the twin prime conjecture.
\end{theorem}

\begin{lemma}[Leibnitz-type identity]
\label{PreobLeibnitz}
Let $L(n)=\log n$ and ${\star}$ be the Dirichlet convolution.
Then
\[
L\times(\mu\star1\star\tau_{\kappa})=(L\mu)\star1\star\tau_{\kappa}+\mu\star(L\tau_{1+\kappa}).
\]
\end{lemma}
\dokvo The lemma follows by differentiating the generating Dirichlet series.

\textbf{2. Proof of Theorem~\ref{Preoblambdacorrel}.}
By Lemma~\ref{PreobLeibnitz} we have
\[
\mu\star(L\tau_{1+\kappa})=L\times(\mu\star1\star\tau_{\kappa})-(L\mu)\star1\star\tau_{\kappa}.
\]
We rewrite the right-hand side using the identities $\mu\star1=1_{n=1}$, $-(L\mu)\star1=\Lambda$ as
\[
L\tau_{\kappa}+\Lambda\star\tau_{\kappa}.
\]
Assuming $n$ squarefree and denoting the number of prime divisors of $n$ by $\omega(n)$
the latter expression is
\[
(\log n)\kappa^{\omega(n)}+(\log n)\kappa^{\omega(n)-1}
\]
with the convention that $0^0=1$. Dividing by $\log N$ we obtain the identity
\[
\frac{\log n}{\log N}\kappa^{\omega(n)}+\frac{\log n}{\log N}\kappa^{\omega(n)-1}
=\sum_{d\mid n}\mu(d)\tau_{1+\kappa}\left(\frac nd\right)\left(\frac{\log n}{\log N}-\frac{\log d}{\log N}\right).
\]
Now the conclusions of the theorem follow from this identity.

\bigskip

\end{document}